\documentclass[11pt]{article}
\usepackage{amsmath}
\usepackage{graphics, amssymb}
\newtheorem{Def}{Definition}[section]
\newtheorem{Thm}[Def]{Theorem}

\newtheorem{Lemma}[Def]{Lemma}

\newtheorem{Corollary}[Def]{Corollary}
\newtheorem{Nota}[Def]{Notation}

\newcommand{\Proof}{{\em{Proof. }}}

\newcommand{\QED}{\ \hfill $\FBox$ \\[1em]}
\newcommand{\spc}{\;\;\;\;\;\;\;\;\;\;}
\newcommand{\R}{\mathbb{R}}

\newcommand{\FBox}{\rule{2mm}{2.25mm}}
\newcommand{\ov}{\overline}
\newcommand{\vol}{\mbox{\rm vol}}
\newcommand{\spec}{\mbox{\rm spec}}

\newcommand{\D}{\mathcal{D}}

\newcommand{\ti}{\tilde}

\newcommand{\ssn}{S_{S^n}}
\newcommand{\1}{\mbox{\rm 1 \hspace{-1.05 em} 1}}
\newcommand{\Tr}{\mbox{\rm{Tr}\/}}

\numberwithin{equation}{section}

\title{A Weighted~$L^2$-Estimate of the Witten Spinor in
Asymptotically Schwarzschild Manifolds}
\author{Felix Finster\thanks{Supported by the Deutsche Forschungsgemeinschaft
within the Priority Program ``Globale Differentialgeometrie''.}
and Margarita Kraus$^*$}
\date{February 2005}

\begin{document}
\maketitle

\begin{abstract}
We derive a weighted $L^2$-estimate of the Witten spinor in
a complete Rie\-mannian spin manifold~$(M^n, g)$ of non-negative scalar curvature
which is asymptotically Schwarz\-schild.
The interior geometry of~$M$ enters this estimate only
via the lowest eigenvalue of the square of the Dirac
operator on a conformal compactification of~$M$.
\end{abstract}

\section{Introduction}
Since Witten's proof of the positive mass theorem~\cite{W, PT, Ba},
spinors have been a valuable tool for the analysis of
asymptotically flat manifolds; see for example~\cite{HGHP,
Herzlich} or the integral estimates
of Riemannian curvature~\cite{BF, FK, FKr}.
These last estimates have the disadvantage that they involve
the isoperimetric constant, which depends on the geometry
in the interior (i.e.\ away from the asymptotic end)
and is therefore in most situations not known. In order to get
curvature estimates which do not involve the isoperimetric
constant, one needs better control of the Witten spinor. This was our
motivation for looking at weighted integral norms of
the Witten spinor.
In order to concentrate on the role of the interior geometry, we
chose the geometry in the asymptotic end as simple as
possible: as the Schwarzschild metric.
Thus the question under consideration is
to which extent the unknown interior geometry
can affect the behavior of the Witten spinor in the asymptotic
end. In this paper, we quantify this effect by an
integral inequality. We find that the effect of
the interior geometry on a suitable
weighted $L^2$-norm is described purely in terms of the
lowest eigenvalue of the square of the Dirac operator
on a conformal compactification of~$M$.

To be more specific, we now describe the problem
and our main results in the most familiar and physically
most interesting case of dimension three.
Thus let~$(M^3, g)$ be a complete Riemannian manifold of
non-negative scalar curvature which is asymptotically
Schwarz\-schild, i.e.\ the metric in the asymptotic
end is
\[ g \;=\; \left(1+\frac{m}{2r} \right)^{4}\:g_0 \:, \]
where~$g_0$ is the Euclidean metric, and~$r=|x|$ is the
Euclidean norm of~$x \in \R^3$.
The Witten spinor is a solution of the
massless Dirac equation which at infinity
goes over to a constant spinor~$\psi_0$ with~$\|\psi_0\|=1$.
More precisely, a Witten spinor~$\psi$ has the following
asymptotics at infinity,
\begin{equation} \label{asy}
\psi(x) \;=\; \left(1+\frac{m}{2r} \right)^{-2} \psi_0 \:+\:
{\mathcal{O}}\!\left(\frac{1}{r^{2}} \right) .
\end{equation}
This asymptotics is used in~\cite{W, PT} for the proof
of the positive mass theorem. In order to get an
estimate of the error term, we point compactify the manifold
with a conformal transformation of the form
\begin{figure}[t]
\begin{center}
\input{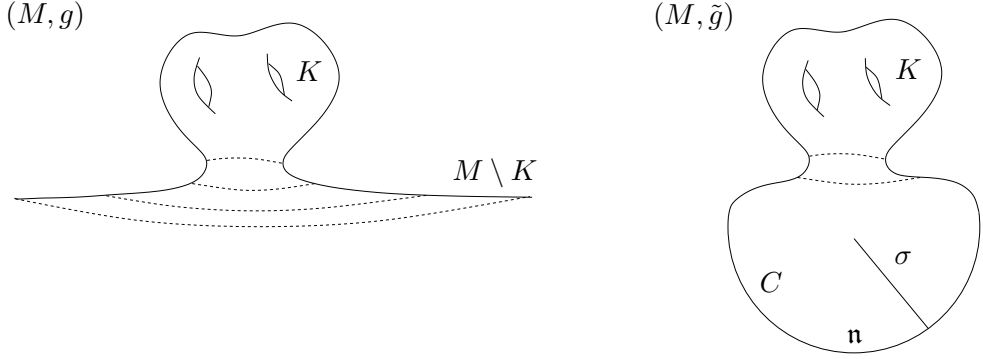}
\caption{The asymptotically Schwarzschild manifold~$(M,g)$ and its
conformal compactification~$(\bar{M}, \tilde{g})$.}
\label{fig1}
\end{center}
\end{figure}
\[ \tilde{g} \;=\; \lambda^2\: g \:, \]
in such a way that the geometry of~$K$ remains unchanged,
the scalar curvature stays non-negative,
and the compactification of the asymptotic end is isometric
to a cap~$C \subset S^n_\sigma$ of a sphere of radius~$\sigma$
(for an illustrating example see Figure~\ref{fig1}).
Then the compactification~$(\bar{M}, \tilde{g})$ is
a closed manifold of non-negative scalar curvature, and
it is clear from the Lichnerowicz-Weitzenb{\"o}ck formula
that the square of the Dirac operator on~$\bar{M}$ is positive,
$\tilde{\D}^2 \geq 0$. Since the scalar curvature is strictly positive
in the spherical cap, we even know that the lowest eigenvalue
is non-zero, $\inf \spec(\tilde{\D}^2)>0$. Moreover,
using methods of spin geometry,
it is possible under various geometric conditions
(for example involving only scalar curvature)
to bound the lowest eigenvalue of~$\tilde{\D}^2$ from below
(see e.g.~\cite{Friedrich, Hi1, Hi2, Lott, Baer, Kraus, Ammann}).
Then the following inequality gives a detailed estimate
of the Witten spinor.

\begin{Thm} \label{thmintro}
There is a constant~$c$ independent of the geometry of~$K$ such that
\[ \int_K \|\psi\|^2 \:d\mu_M \:+\: \int_{M \setminus K} \left\| \psi - \left(1+\frac{m}{2r} \right)^{-2} \psi_0
\right\|^2\:\lambda\:d\mu_M \;\leq\; \frac{c}{\inf \spec({\tilde{\D}}^2)} \: . \]
\end{Thm}
In the course of proving this theorem, we derive an identity
involving Witten spinors, which is of some interest in its own,
as we now explain. We choose an orthonormal basis~$(\psi_{0,i})_{i=1,\ldots,4}$
of the spinors at infinity and consider the corresponding family~$\psi_i$
of Witten spinors,
\begin{equation} \label{fWitten}
\D\psi_i = 0 \:,\qquad \lim_{|x| \to \infty} \psi_i(x) = \psi_{0,i} \spc (i=1,\ldots,4).
\end{equation}
We let~$G$ and~$G_{S^3_\sigma}$ be the Green's functions of the square
of the Dirac operator on~$\bar{M}$ and the sphere~$S^3_\sigma$,
respectively, and denote their integral kernels by~$G(x,y)$
and~$G_{S^3_\sigma}(x,y)$. The next theorem expresses the weighted
$L^2$-norms of the Witten spinors~(\ref{fWitten}) in terms of
the difference of these integral kernels, with an explicit
error term.
\begin{Thm} \label{thmintro2}
The Witten spinors satisfy for sufficiently large~$R$
the identity
\begin{eqnarray*}
\lefteqn{ \int_K \sum_{i=1}^4 \|\psi_i\|^2 \:d\mu_M \:+\: \int_{M \setminus K}
\sum_{i=1}^4 \left\| \psi_i - \left(1+\frac{m}{2r} \right)^{-2} \psi_{0, i}
\right\|^2\:\lambda\:d\mu_M } \\
&=& 64 \pi^2 \sigma^4\: \lim_{\mathfrak{n} \neq y \to \mathfrak{n}}
\Tr (G(\mathfrak{n}, y) - G_{S^3_\sigma}(\mathfrak{n}, y)) \\
&&+ 4 \int_{B_R(0)} \frac{2 \sigma^2}{\sigma^2+r^2}\: d^3x
\:-\:4 \int_{B_R(0) \cap (M \setminus K)}
\left( 1 + \frac{m}{2r} \right)^2 \:\lambda\: d^3x\:,
\end{eqnarray*}
where~$d^3x$ is the Lebesgue measure on~$\R^3$. Here
we assume that the asymptotic end~$M \setminus K$ is
diffeomorphic to~$\R^3 \setminus B_\rho(0)$ and work in the corresponding
chart.
\end{Thm}
This {\em{equality}} clearly gives finer
information than the inequality of Theorem~\ref{thmintro}.
The interesting point is that the interior geometry enters
only via the Green's function~$G$. We learn that
the influence of the interior geometry
on the weighted~$L^2$-norm is described precisely by the
behavior of~$G(\mathfrak{n},y)$ as~$y \to \mathfrak{n}$.

In this paper, we prove the analog of Theorem~\ref{thmintro}
and Theorem~\ref{thmintro2} in general dimension.
For the proof we work as in~\cite{FK} with the spinor operator,
which is composed of a basis of Witten spinors~(\ref{fWitten}).
Our first step is to get a connection between the
conformally transformed spinor operator and
a quadratic expression in the Dirac Green's
function on~$\bar{M}$ (see Section~\ref{sec3}).
In Section~\ref{sec4} we find that, after subtracting
suitable counter terms, we can integrate this expression
over~$\bar{M}$ to obtain the Green's function~$G$
of the square of the Dirac operator minus suitable counter
terms. Then our task becomes to analyze the behavior of~$G$
near the pole~$\mathfrak{n}$ of the spherical cap.
This is done in Section~\ref{sec5}, where we estimate the
difference of~$G$ and the corresponding Green's function
on the sphere using Sobolev techniques. In
Section~\ref{secGFS}, we compute the Green's functions
on the sphere explicitly. Finally, in Section~\ref{sec7} we
combine the results of Sections~\ref{sec3} and~\ref{sec4}
to obtain an identity for an integral of the trace of the
spinor operator (Theorem~\ref{Thmw1}).
Using a positivity argument together with the estimates
of Sections~\ref{sec5} and~\ref{secGFS}, we then
conclude our main results (Corollary~\ref{corend},
Theorem~\ref{thmspop} and Corollary~\ref{thmwitten}).

We finally remark that our methods work similarly
also for harmonically flat manifolds (see~\cite{SY2} for the definition),
except for Lemma~\ref{partialwave}, where an angular dependence
of the function~$\lambda$ would lead to a ``mixing'' of the
angular momentum modes, making the situation more complicated.
Since we did not find a simple argument to overcome this problem,
we here restrict attention to a Schwarzschild end.
A generalization to harmonically flat manifolds would be
desirable in view of the fact that the metric of every asymptotically
flat manifold can be made harmonically flat by
an arbitrarily small perturbation~\cite{SY2}.

\section{The Conformal Compactification} \label{seccc}
Let~$(M^n,g)$ be a complete Riemannian manifold of dimension~$n\geq3$
with non-negative scalar curvature.
For simplicity, we assume that the manifold has one asymptotic end which
is isometric to Schwarzschild.
By rescaling, we can assume without loss of generality
that the ADM mass is equal to two.
\begin{Def} \label{asymp Schwarzschild}
A complete Riemannian manifold $(M^n,g)$ of dimension $n\geq3$ is said to be
{\bf{asymptotically Schwarzschild}} if there is a parameter $\rho > 0$, a
compact set $K \subset M$ and a diffeomorphism $$\phi : M \setminus K \to \R^n \setminus \overline{B_{\rho}(0)} $$
such that
\begin{equation} \label{Schmet}
\phi_* g \;=\; \left(1+\frac{1}{|x|^{n-2}}\right)^{\frac{4}{n-2}}\:g_0\:.
\end{equation}
Here~$B_r(0)$ denotes an open Euclidean ball in~$\R^n$, and
$g_0$ is the Euclidean metric on~$\R^n$.
\end{Def}
In the asymptotic end it is most convenient to work in the chart~$(\phi,\,
M \setminus K)$; we use the notation
\[ r(x) \;=\; |\phi(x)| \quad {\mbox{for $x \in M \setminus K$}}. \]
The metric~(\ref{Schmet}) is obviously conformally
equivalent to the Euclidean metric. Moreover, the Euclidean~$\R^n$
is conformal to a sphere~$S^n_\sigma$ of radius~$\sigma$ with the
north pole removed.
This is is seen explicitly in the usual chart obtained
by stereographic from the north pole, where
\begin{equation} \label{gSn}
g_{S^n_\sigma} \;=\;
\left( \frac{2 \sigma^2}{\sigma^2+r^2} \right)^2\: g_0 \:.
\end{equation}
Therefore, we can arrange by a conformal transformation that the
asymptotic end is isometric to the cap of a sphere of radius~$\sigma$
with the north pole removed. Furthermore, we want to keep the metric
inside~$K$ unchanged, and we want to preserve the positivity of scalar
curvature. In the next lemma we construct a function~$\lambda$ such that
the conformal transformation
\begin{equation} \label{cc}
\tilde{g} \;=\; \lambda^2\: g
\end{equation}
gives us a metric with the desired properties.

\begin{Lemma} \label{lemmacap}
There is a function~$\lambda \in C^\infty(M)$ satisfying the
following conditions:
\begin{description}
\item[(i)] $\displaystyle \lambda_{|K} \;\equiv\; 1$
\item[(ii)] For some radius $R >\rho$,
\begin{equation} \label{ldef}
\lambda(x) \;=\; \left(\frac{2 \sigma^2}{\sigma^2+r(x)^2}\right)\cdot
 \left(1+\frac{1}{r(x)^{n-2}}\right)^{-\frac{2}{n-2}}
\qquad {\mbox{on~$\phi^{-1}(\R^n \setminus B_{R}(0))$}}\:.
\end{equation}
\item[(iii)] The scalar curvature corresponding to the
metric~(\ref{cc}) is non-negative.
\end{description}
More specifically, we can arrange that
\begin{equation} \label{scaling}
\rho \;\leq\; \sigma \;\leq \; R \;\leq\; c(n)\: (\rho+1)
\end{equation}
with a constant~$c$ which depends only on the dimension.
\end{Lemma}
{\Proof} It is more convenient to write the metric in the
asymptotic end as~$\tilde{g}= \mu^2 g_0$ with
\[ \mu(x) \;=\; \left(1+\frac{1}{r(x)^{n-2}}\right)^{\frac{2}{n-2}}\:
\lambda(x) \:. \]
Our first attempt is to define~$\mu$ piecewise. To this end, we
choose~$R_* = \rho + C(n)$ with a suitable constant~$C(n)$ and set
\[ \mu \;=\; \mu(r(x)) \;=\; \left\{ \begin{array}{cl}
\displaystyle \left(1+\frac{1}{r^{n-2}}\right)^{-\frac{2}{n-2}} & {\mbox{if $r<R_*$}} \\[.7em]
\displaystyle \frac{2 \sigma^2}{\sigma^2+r^2} & {\mbox{if $r \geq R_*$}} 
\;. \end{array} \right. \]
In order to make this function continuous at~$r=R_*$, we let
\[ \sigma \;=\; R_* \left[ 2 \left(1+\frac{1}{R_*^{n-2}}
\right)^{-\frac{2}{n-2}} - 1 \right]^{-\frac{1}{2}} \:. \]
By choosing~$C(n)$ sufficiently large, we can arrange that
the square bracket is uniformly bounded from above and below.
In the region~$r<R_*$, the metric~$\tilde{g}$ coincides
with~$g$, whereas in the
region~$r>R_*$, $\tilde{g}$ is the metric
on~$S^n_\sigma$. Obviously, in both of these regions the conformally
transformed metric has all the required properties.
Unfortunately, the function~$\mu$ is not smooth at~$r=R_*$.
More precisely, a short calculation shows that the first derivative
of~$\mu$ makes a negative jump, i.e.
\[ \lim_{R_* < r \to R_*} \mu'(r)
- \lim_{R_* > r \to R_*} \mu'(r) \;<\; 0\:. \]
As a consequence, the scalar curvature~$\tilde{s}$ corresponding
to~$\tilde{g}$, given by the formula
\[ \tilde{s} \;=\; 4\: \frac{n-1}{n-2}\: \mu^{-\frac{n+2}{2}} \:
\Delta \mu^{\frac{n-2}{2}} \]
(where~$\Delta = - \nabla^i \nabla_i$ is the Laplacian in~$\R^n$)
is positive at~$r=R_*$ in the distributional sense.
By mollifying~$\mu(r)$ in a small neighborhood of~$r=R_*$
we can thus arrange that~$\tilde{s} \geq 0$.
We finally set~$R=\sigma+1$.
\QED

For clarity, we denote the manifold~$M$ with metric~(\ref{cc})
by~$(\tilde{M}, \tilde{g})$. Then by~{\bf{(i)}}, the manifolds~$M$ and $\tilde{M}$ are isometric on~$K$.
By~{\bf{(ii)}}, $\tilde{M}$ is on~$\phi^{-1}(\R^n \setminus B_{R}(0))$
isometric to the cap of the sphere~$S^n_\sigma$ with the north
pole~$\mathfrak{n}$ removed.
We denote the geodesic distance from the north pole by $d$,
\begin{equation} \label{distdef}
d \::\: S^n_\sigma \to \R_+
\end{equation}
and let~$B_s(\mathfrak{n})$ be the geodesic balls of radius~$s$
around the north pole. Then
\[ \phi^{-1}(\R^n \setminus B_{R}(0))
\;\cong\; \left( B_\delta(\mathfrak{n}) \setminus \{\mathfrak{n}\}
\right) \subset S^n_\sigma
\quad {\mbox{with}} \quad \delta = 2 \sigma \: \arctan \left(
\frac{\sigma}{R} \right) . \]

We now compactify~$\tilde{M}$ by adding the north pole.
The resulting manifold, denoted by~$(\ov{M},\tilde{g})$,
is called the {\em{conformal compactification}} of~$(M, g)$.
For~$r \geq R$ we set
\[ C_r \;=\; \ov{(\varphi^{-1}(\R^n \setminus B_r(0)),
\tilde{g})} \subset \bar M \:. \]
We also identify~$C_{r} $ via the isometry of the stereographic projection
with a closed subset of~$S^n_\sigma$.
We refer to the set~$C \equiv C_{R}$ as the {\em{spherical cap}} of $(\ov{M},\tilde{g})$. We always identify it with the
set~$B_\delta(\mathfrak{n}) \subset S^n_\sigma$.

We remark that there are also conformal compactifications
with the above properties~{\bf{(i)}}--{\bf{(iii)}},
for which~$\sigma$ is arbitarily large, thus violating~(\ref{scaling}).
However, such conformal compactifications do not seem to give good
estimates of the Witten spinor, and we shall not consider them here.

\section{The Spinor Operator and the Dirac Green's Function} \label{sec3}
From now on we need to assume that~$(\bar{M}, \tilde{g})$ is a {\em{spin}} manifold. This assumption is no restriction in dimension three,
whereas in general it poses a constraint for the topology of~$\bar{M}$.
As a consequence, the manifold~$M$ is also spin. In fact,
taking out the point~$\mathfrak{n}$ and performing a conformal
transformation, every spin structure on~$\bar{M}$
induces a spin structure on~$M$.
We fix corresponding spin structures on~$M$ and~$\bar{M}$ throughout.

We let $\Sigma$ and $\tilde{\Sigma}$ be conformally equivalent spinor bundles
over~$(M,g)$ and $(M,\tilde{g})$ and denote the corresponding Dirac operators
by~$\D$ and~$\tilde{\D}$.  According to \cite{Hitchin74, Hi1} there is a
fiberwise isometry $\Sigma \to \tilde{\Sigma}, \,
\psi \mapsto \tilde{\psi}$ such that
\begin{equation} \label{ctrans}
\tilde{\D} \tilde{\psi}
\;=\; \lambda^{-\frac{n+1}{2}}\,\widetilde{ \left( \D (\lambda^{\frac{n-1}{2}}
\psi ) \right) }.
\end{equation}
In the coordinates induced by the diffeomorphism $\phi$ of
Definition~\ref{asymp Schwarzschild}, we choose a family of
constant spinors $\psi_{0,i}$ such that
$\psi_{0,1},\ldots,\psi_{0,N}$, $N=2^{[n/2]}$, is an orthonormal basis
in the asymptotic end~$M \setminus \varphi^{-1}(\R^n \setminus B_R(0))$
and consider the boundary value problem
\begin{equation}\label{bv}
 \D\psi_i = 0 \:,\spc \lim_{|x| \to \infty} \psi_i(x) = \psi_{0,i}\:.
\end{equation}
This boundary value problem was first considered in~\cite{W},
its solutions are called {\em{Witten spinors}}. The existence and
uniqueness of the Witten spinors was proven in~\cite{PT, Ba}.
They decay at infinity as
\[ \psi_i \;=\; \psi_{0,i}+\mathcal{O}(r^{2-n}), \quad \partial_j \psi_i \;=\; \mathcal{O}(r^{1-n}),
\quad \partial_k\partial_l\psi_i=\mathcal{O}(r^{-n})\:. \]
In~\cite{FK} the {\em{spinor operator}}~$\Pi_x$ was introduced,
which we now slightly generalize, using the following

\begin{Nota}
 Let $E_1\to X, \, E_2 \to X$ be vector bundles
 over the manifold $X$. By
 $$E_1 \boxtimes E_2 \to X \times X$$
 we denote the vector bundle
 $$\pi_1^* E_1 \otimes \pi^*_2 E_2$$
 with the projections $\pi_i : X \times X \to X$ to the
 $i^{\mbox{\scriptsize{th}}}$ factor.
\end{Nota}

\begin{Def} \label{defso}
 Let $\psi_1, \ldots, \psi_N$ be a family of solutions of (\ref{bv}).
Then the {\bf{spinor operator}}~$\Pi \in \Gamma(\Sigma \boxtimes \Sigma^*)$
is defined by
$$ \Pi(x,y)\;=\; \sum_{i=1}^N \langle\psi_i(y),\;.\;\rangle \:\psi_i(x)$$
\end{Def}
This definition reduces to the spinor operator as used in \cite{FK} if
$x$ and $y$ are equal,
\[ \Pi(x) \;:= \;\Pi (x,x) \;:\; \Sigma_x \to \Sigma_x \:. \]
From the boundary values (\ref{bv}) it is obvious that
 $$\lim_{|x| \to \infty} \Pi(x) \;=\; \1 \:.$$
When considering conformal transformations of the spinor operator,
we must keep in mind that the transformed Witten spinor should
again be a solution of the Dirac equation.
According to~(\ref{ctrans}), this leads us to the following definition.
\begin{Def} \label{defctso}
On the manifold $(M,\tilde{g})$ with $\tilde{g}=\lambda^2g $ we
define the {\bf{conformally transformed spinor
operator}}~$\tilde{\Pi} \in \Gamma(\tilde{\Sigma} \boxtimes
 \tilde{\Sigma}^*)$ by
$$\tilde{\Pi}(x,y) \;:=\; \lambda^{\frac{1-n}{2}}(x)\: \lambda^{\frac{1-n}{2}}(y)
 \sum_{i=1}^{N} \langle\tilde{\psi_i}(y),\;.\;\rangle \:
\tilde{\psi_i}(x) \:, $$
where $\psi_1, \ldots, \psi_N$ are the solutions of (\ref{bv}).
\end{Def}
In Euclidean space, the Witten spinors are constant,
and therefore~$\Pi(x) \equiv \1$. Applying the above definition
to the metric~(\ref{asymp Schwarzschild}),
we obtain for the spinor operator~$\Pi_{\mbox{\scriptsize{Sch}}}$
in Schwarzschild the explicit expression
\begin{equation} \label{PiSch}
\Pi_{\mbox{\scriptsize{Sch}}}(x) \;=\; \left(1 +
\frac{1}{r(x)^{n-2}} \right)^{-2\: \frac{n-1}{n-2}}\: \1_{\Sigma_x}\:.
\end{equation}

In the remainder of this section we will establish a
connection between the conformally transformed spinor
operator on~$(\bar{M}, \tilde{g})$
and the Green's function of the Dirac operator.
\begin{Def}
Let $X$ be a spin manifold and~$\Delta$ the diagonal of~$X \times X$,
\[ \Delta \;=\;  \{ (x,x) {\mbox{ with~$x \in M$}} \} \subset X \times X\:. \]
We let~$S_X$ be a smooth section in the bundle $\Sigma \boxtimes \Sigma^* | ((X \times X) \setminus \Delta)$,
$$ S_X :\; X \times X \setminus \Delta \longrightarrow \Sigma_X \boxtimes \Sigma^*_X, $$
and also consider $S_X$ as the integral kernel of a corresponding operator
acting on the compactly supported, smooth spinors on $X$ by
$$ (S_X \psi ) (x) = \int_{X \setminus \{x\}} S_X (x,y)\:\psi(y) dy\: . $$
$S_X$ is called the {\bf{Green's function of the Dirac operator}} $\D$ on X if
it satisfies the distributional equation
$$\D_{X,x} \, S_X(x,y) \;=\; \delta(x,y)\:. $$
\end{Def}
This distributional equation can be stated equivalently by
the condition that
$$\int_X S_X(x,y)\:\D_x \psi(x) \:d^nx\;=\; \psi(y)$$
for all compactly supported, smooth sections in the spinor bundle.
As is easily verified by a direct computation,
the Green's function on~$\R^n$ with the Euclidean metric is given by
\begin{equation} \label{SEuclid}
S_{\R^n}(x,y) \;=\; -\frac{1}{\omega_{n-1}}\frac{x-y}{|x-y|^n} \:,
\end{equation}
where $\omega_{n-1}$ is the volume of $S^{n-1}$.
The Green's functions of conformally flat spaces can be computed
using the following transformation law for the Green's function under
conformal changes.
\begin{Lemma}\label{conformal change}
Let~$(X,g)$ and~$(\tilde{X}, \tilde{g})$ be two manifolds with 
conformally equivalent metrics, $\tilde{g} = \lambda^2 g$.
Then the corresponding Green's functions~$S$ and~$\tilde{S}$
are related by
$$ S_{\tilde{X}}(x,y) \;=\;
\lambda^{\frac{1-n}{2}}(x) \: \lambda^{\frac{1-n}{2}}(y) \; S_X(x,y)\:. $$
\end{Lemma}
\Proof Let~$S$ be the Green's function of the Dirac operator
$(X,g)$, i.e.
\[ \D_x \!\int_x S_X(x,y) \, \psi(y)\:dy \;=\; \psi(x) \:. \]
Then from~(\ref{ctrans}),
$$ \tilde{\D}_x \, \lambda^{-\frac{n-1}{2}}(x)
\left( \int_{\tilde{X}} S_X(x,\tilde{y}) \, \psi(\tilde{y}) \,
\lambda^{-n}(\tilde{y}) \: d\tilde{y} \right) \tilde{} \;=\;
\lambda^{-\frac{n+1}{2}}(x) \, \tilde{\psi}(x)$$
and thus
$$ \tilde{\D}_x \left( \int_{\tilde{X}} \lambda^{-\frac{n-1}{2}}(x) \,
S_X(x,\tilde{y}) \, \lambda^{-\frac{n-1}{2}}(\tilde{y}) \,
\lambda^{-\frac{n+1}{2}}(\tilde{y}) \, \psi(\tilde{y}) \:
d\tilde{y} \right) \tilde{}  \;=\; \lambda^{-\frac{n+1}{2}}(x) \,
\tilde{\psi}(x)\:, $$
where  $\tilde{\D}$ and $d\tilde{y}$ denote the
Dirac operator and the volume element on $\tilde{X}$, respectively. \QED

\begin{Thm}\label{1.Theorem}
Let $(M,g)$ be an asymptotically Schwarzschild manifold and
$(\bar{M},\tilde{g})$ its conformal compactification. Then
the conformally transformed spinor operator $\tilde{\Pi}$ and
the Green's function $\tilde{S}_{\bar{M}}$ of the Dirac operator on
$(\bar{M},\tilde{g})$ satisfy the following identity,
$$ \tilde{\Pi}(x,y) \;=\; \omega^2_{n-1}\:(2 \sigma^2)^{n-1}
\:\tilde{S}_{\bar{M}}(x,\mathfrak{n}) \:
\tilde{S}_{\bar{M}}(\mathfrak{n},y) \:, $$
where~$\mathfrak{n} \in C \subset (\bar{M},\tilde{g})$ is the
north pole in the spherical cap of the compactification.
\end{Thm}
\Proof Let~$B_{\varepsilon}(\mathfrak{n})\subset C$ be a ball of
radius $\varepsilon$ around $\mathfrak{n}$, $\varepsilon<\delta$. Then
\begin{eqnarray*}
 \int_{\bar{M}\setminus B_{\varepsilon}(\mathfrak{n})}
 \tilde{\Pi}(x,y)  \left( \tilde{\D}\,\tilde\psi
 \right)(y) \:dy &=& \sum_{i=1}^N \; \int_{\bar{M}\setminus B_{\varepsilon}(\mathfrak{n})}
 \textrm{div} \:V_{\tilde{\psi},\psi_i}(x) \:dx\; \hat{\psi_i}(y) \\
 &=& \sum_{i=1}^N \int_{\partial
 B_{\varepsilon}(\mathfrak{n})} \langle \tilde{\psi}, n_{\varepsilon} \cdot \hat{\psi_i} \rangle(x) \:dx\;\hat{\psi}_i(y) \:,
\end{eqnarray*}
where~$\hat{\psi}_i = \lambda^{-\frac{(n-1)}{2}}\cdot \tilde{\psi}_i$,
and ~$n_{\varepsilon}$ is the outer normal on $\partial B_{\varepsilon}(\mathfrak{n})$.
Here the vector field~$V_{\psi,\psi_i}$ is defined by
$$ \tilde{g}(V_{\psi,\psi_i},w) \;=\; \langle\tilde{\psi},w \cdot\tilde{\psi_i}\rangle \, \lambda^{-\frac{(n-1)}{2}}(x) \:. $$
Similarly, we obtain
\begin{eqnarray}\label{2stars}
 && \int_{\bar{M}\setminus B_{\varepsilon}(\mathfrak{n})} \int_{\bar{M}\setminus
 B_{\varepsilon'}(\mathfrak{n})} \left\langle \tilde{\D} \tilde{\varphi_1}(x)
, \tilde{\Pi}(x,y)\: \tilde{\D} \tilde{\varphi_2}(y) \right\rangle  \:dx \:dy\\
 && = \;
 \sum_{i=1}^N \left( \; \int_{\partial B_{\varepsilon}(\mathfrak{n})} \left\langle
 \tilde{\varphi_1}, n_{\varepsilon} \cdot \hat{\psi}_i \right\rangle\!(x) \:dx \right)
 \cdot \left( \int_{\partial B_{\varepsilon'}(\mathfrak{n})}
 \left\langle n_{\varepsilon'} \cdot \hat{\psi}_i, \tilde{\varphi_2} \right\rangle\!(y)\:dy
 \right) .
\end{eqnarray}
In order to investigate the limit $\varepsilon \to 0$,
we consider the trivialization of the spinor bundle in a
neighborhood of the north pole given by stereographic projection
from the south pole $-\mathfrak{n}$.
The change of the charts of $S^n_\sigma \setminus
(\mathfrak{n} \cup -\mathfrak{n})$ given by the stereographic
projections from the north and south pole is given by
$\omega(x) = \sigma^2 \frac{x}{|x|^2}$ with
differential $d\omega_x (v) = \frac{\sigma^2}{|x|^2} S(\frac{x}{|x|})v$,
where $S(\frac{x}{|x|})$ denotes the reflection at $(\R \cdot
x)^{\perp}$. Thus, for an orthonormal basis $(e_1,\ldots,e_n)$ of
$\R^n$, the vector fields $\tilde{e}_i : \R^n \setminus 0 \to \R^n$, $x \mapsto
\frac{\sigma^2+|x|^2}{2 \sigma^2}S(\frac{x}{|x|})\: e_i$ are extendable to the sphere
$S^n_\sigma \setminus (-\mathfrak{n})$ as orthonormal vector fields.  The lift of
$S(\frac{x}{|x|}) \in \mathrm{O}(n)$ to $\mathrm{Pin}(n)$ is
given by Clifford multiplication with
$\frac{x}{|x|}= \frac{2 \sigma^2}{\sigma^2+|x|^2} d \rho_{\rho^{-1}(x)} (n_\varepsilon)$
for $x \in \rho (B_\varepsilon(\mathfrak{n}))$, and therefore the Clifford products
$(n_{\varepsilon} \cdot \tilde{\psi}_i)_{i=1,\ldots,N}$ are extendable to
the north pole as orthonormal basis $\tilde{\psi}_{0,1}(\mathfrak{n}), \ldots ,
\tilde{\psi}_{0,2^{[\frac n2]}}(\mathfrak{n})$.

From the asymptotic behavior of the solutions of (\ref{bv}) it follows that
$$n_{\varepsilon} (x) \cdot \hat{\psi}_i(x) \;=\; (2\sigma^2)^{\frac{1- n}{2}} r'(x)^{1-n}\cdot
\tilde{\psi}_{0,i}(\mathfrak{n}) + \mathcal{O}(|r'(x)|^{2-n}) \quad \textrm{
with } \quad
\langle \tilde{\psi}_{0,i}(\mathfrak{n}),\tilde{\psi}_{0,j}(\mathfrak{n}) \rangle
\;=\; \delta_{ij} \:, $$
where $(r', \Omega) : S^n_\sigma \to \R^n$ denotes the coordinates of the stereographic
projection from the south pole. Therefore,
\begin{eqnarray*}
\lim_{\varepsilon \to 0} \int_{\partial
 B_{\varepsilon}(\mathfrak{n})} \langle \tilde{\varphi}(x), n_{\varepsilon}
 \cdot \hat{\psi}_i(x) \rangle\:  dx
&=& \lim_{\varepsilon \to 0} \varepsilon^{1-n} (2 \sigma^2)^{\frac{n-1}{2}}
 \vol(S^{n-1}_{\varepsilon}) \:\langle\tilde{\varphi}(\mathfrak{n}),  \tilde{\psi}_{0,i}(\mathfrak{n}) \rangle 
\:+\: \mathcal{O}(\varepsilon)\\
&=& \omega_{n-1} (2 \sigma^2)^{\frac{n-1}{2}}\:
\langle \tilde{\varphi}(\mathfrak{n}),\tilde{\psi}_{0,i}(\mathfrak{n}) \rangle \:.
\end{eqnarray*}
Now we can in (\ref{2stars}) take the limit~$\varepsilon \to 0$ to obtain
\begin{eqnarray*}
\lefteqn{ \int_{\bar{M} \times \bar{M}} \left\langle 
\tilde{\D} \tilde{\varphi_1}(x), \tilde{\Pi}(x,y)\: \tilde{\D} \tilde{\varphi_2}(y)
\right\rangle\: dx \: dy \;=\; (2 \sigma^2)^{n-1}
 \left\langle\tilde{\varphi_1}(\mathfrak{n}),\tilde{\varphi_2}(\mathfrak{n})
 \right\rangle \cdot \omega_{n-1}^2 } \\
&& = (2 \sigma^2)^{n-1} \int_{\bar{M} \times \bar{M}}
\left\langle S_{\bar{M}}(\mathfrak{n},x)\, \tilde{\D} \, \varphi_1(x),
 S_{\bar{M}}(\mathfrak{n},y)\, \tilde{\D} \,\varphi_2(y) \right\rangle
\:dx\:dy  \cdot \omega_{n-1}^2 \:.
\end{eqnarray*}

\vspace*{-2.05em} \QED 

\section{The Green's Function of the Square of the Dirac Operator} \label{sec4}
Let us outline our strategy. Our goal is to derive
weighted $L^2$-estimates of the Witten spinors. Since the
spinor operator is composed of the Witten spinors (see
Definition~\ref{defctso}), the expression
\[ \int_{\bar{M}} \: \Tr\: \tilde{\Pi}(x) \: dx \]
is of interest, where ``Tr'' denotes the trace on
the $N$-dimensional vector space $\Sigma_x$ (for details
see Section~\ref{sec7}). Using the formula from Theorem~\ref{1.Theorem}
and the cyclicity of the trace, we are led to the integral
\begin{equation} \label{int1}
\int_{\bar{M}} \: \tilde{S}_{\bar{M}}(\mathfrak{n},x)\: \tilde{S}_{\bar{M}}(x,\mathfrak{n})
\: dx \:.
\end{equation}
If the last argument of the second factor~$\tilde{S}_{\bar{M}}$ were
different from~$\mathfrak{n}$, we could immediately carry out the integral,
\begin{equation} \label{int2}
\int_{\bar{M}} \: \tilde{S}_{\bar{M}}(\mathfrak{n},x)\: \tilde{S}_{\bar{M}}(x,y)
\: dx \;=\; G(\mathfrak{n},y)\spc (y \neq \mathfrak{n})\:,
\end{equation}
where~$G$ denotes the Green's function of the Dirac operator squared,
\begin{equation} \label{Gdef}
\tilde{\D}^2_x G(x,y) = \delta(x-y)\:.
\end{equation}
This simple argument suffers from the problem that the integrals
in~(\ref{int1}) and~(\ref{int2}) diverge. Namely, since the
order of the pole of~$\tilde{S}_{\bar{M}}$ is expected to be the same
as in Euclidean space~(\ref{SEuclid}), we find that the product
of the Green's functions should have a non-integrable pole of the form
\[ \tilde{S}_{\bar{M}}(\mathfrak{n},x)\: \tilde{S}_{\bar{M}}(x,\mathfrak{n})
\;\sim\; \frac{1}{d^{2n-2}}\:, \]
where~$d$ denotes the geodesic distance from the north
pole~(\ref{distdef}).
Despite this problem, one can hope that the
above argument works if suitable
functions are subtracted from the integrands in order compensate the
singularities. This leads us to conjecture a relation of the following form,
\[ \int_{\bar{M}} \left(\tilde{S}_{\bar{M}}(\mathfrak{n},x)\: \tilde{S}_{\bar{M}}(x,\mathfrak{n})
- {\mbox{(counter terms)}} \right) dx
\;=\; \lim_{\mathfrak{n} \neq y \to \mathfrak{n}} \left(G(\mathfrak{n},y) - {\mbox{(counter terms)}} \right) . \]
In order to specify the counter terms, we let
$$ \chi_\delta \;:=\; \chi_{C}$$
be the characteristic function of the spherical cap. Now we take the Dirac Green's function of the sphere
and multiply it by~$\chi_\delta$, so that it is supported inside
the cap~$C$ around the north pole. Using the isometry with
the spherical cap of~$\bar{M}$, we can lift~$S_\delta$ to~$\bar{M}$,
\[ S_\delta \::\: \bar{M} \times \bar{M} \setminus \Delta \;\to\;
\Sigma_{\bar{M}}\,\boxtimes \,\Sigma_{\bar{M}}^* \quad{\mbox{with}}\quad
S_\delta(x,y) \;=\; \chi_\delta(x) \:S_{S^n_\sigma}(x,y)\: \chi_\delta(y) \:. \]
Finally, we set
\begin{equation} \label{Gddef}
G_\delta(y) \;=\; \int_{\bar{M}} S_\delta(\mathfrak{n},x)\:
S_\delta(x,y)\: dx\:.
\end{equation}

\begin{Thm}\label{theorem1}
The Dirac Green's function $\tilde{S}$ on the compact manifold
$(\bar{M},\tilde{g})$ satisfies the relation
\[ \int_{\tilde{M}}\left(\tilde{S}_{\bar{M}}(\mathfrak{n},x)\:
\tilde{S}_{\bar{M}}(x,\mathfrak{n})
- S_\delta(\mathfrak{n},x)\: S_\delta(x,\mathfrak{n}) \right) dx \;=\;
\lim_{\mathfrak{n} \neq y\to \mathfrak{n}}
\left[G(\mathfrak{n},y)- G_\delta (y) \right] . \]
\end{Thm}
In Section~\ref{secGFS} the counter terms are computed more explicitly, see
Lemma~\ref{lemma3}.

The remainder of this section is devoted to the proof of the above
theorem. For any $y\in B_{\delta/4}(\mathfrak{n})$, we introduce the function
$f_y:\bar{M}\to\mathbb{R}$ by
\begin{equation}\label{a}
f_y(x)=\tilde{S}_{\bar{M}}(x,y)-{S}_\delta(x,y)\:.
\end{equation}

\begin{Lemma}\label{lemma1}
There is a constant~$c$ such that
\begin{eqnarray*}
\|f_y\|_\infty &<& c \quad\;\;\;\qquad \mbox{\rm for all}\,\,y\in
B_{\frac{\delta}{4}}(\mathfrak{n}) \\
\lim_{y\to \mathfrak{n}} f_y(x) &=& f_{\mathfrak{n}}(x)
\qquad {\mbox{\rm{uniformly in~$x \in \bar{M}$}}}.
\end{eqnarray*}
\end{Lemma}
\Proof We let~$\mu \in C^\infty_0(\R)$ be a non-negative
test function~$\mu \in C^\infty_0(\R)$ with
$\mu_{|[0,\frac{1}{2}]} \equiv 1$, $0 \leq \mu \leq 1$
and ${\mbox{supp}}\, \mu \subset (-1,1)$
and define the function~$\eta : S^n_\sigma \to \R$ by
\begin{equation} \label{edef}
\eta(x) = \mu\left(\frac{d(x)}{\delta} \right) ,
\end{equation}
where~$d$ is again the geodesic distance from
the north pole~(\ref{distdef}). This function is supported inside
the spherical cap and is identically equal to one in a
the neighborhood~$B_{\delta/2}(\mathfrak{n})$ of the north pole.
For~$y \in B_{\delta/4}(\mathfrak{n})$ we set
\begin{equation*}
g_y(x) \;=\; \tilde{S}_{\bar{M}}(x,y)
- \eta(x)\: {S}_\delta(x,y)\:.
\end{equation*}
Then
\begin{equation*}
\tilde{\D} \,g_y(x) = -(\mbox{grad}\,\eta(x))\: S_{S^n_\sigma}(x,y)\:.
\end{equation*}
Since the function~$\mbox{grad}\,\eta$ is supported in
the annulus~$B_\delta(\mathfrak{n}) \setminus B_{\delta/2}(\mathfrak{n})$,
whereas~$y \in B_{\delta/4}(\mathfrak{n})$, and using that the
Green's function on the sphere is smooth away from the pole~$x=y$
(see Section~\ref{secGFS} for details), we conclude that
$\tilde{\D}\, g_y \in C^\infty(\tilde{M})$. Applying standard
elliptic regularity theory, we obtain that $g_y \in C^\infty(\bar{M})$ and
that it is uniformly bounded in~$y$. Again using 
that~$S_{S^n_\sigma}$ is smooth away from the diagonal,
we conclude that $\tilde{S}_\delta(x,y)$ is
smooth for $x\in B_\delta(\mathfrak{n})\setminus B_{\delta/2}(\mathfrak{n})$, and is bounded
uniformly in $x$ and $y$. We finally note that
\begin{equation*}
f_y(x)=g_y(x)+(1-\eta_\varepsilon (x))\tilde{S}(x,y)\: .
\end{equation*}

\vspace*{-1.85em}\QED

Using~(\ref{a}), we decompose the product of Green's functions
as follows,
\begin{eqnarray}\label{h}
\lefteqn{ \tilde{S}_{\bar{M}}(\mathfrak{n},x)\tilde{S}_{\bar{M}}(x,y)
\;=\; \tilde{S}_{\bar{M}}(x,\mathfrak{n})^\ast \:\tilde{S}_{\bar{M}}(x,y) }
\nonumber \\
&=&({S}_\delta(x,\mathfrak{n})+f_{\mathfrak{n}}(x))^\ast\: ({S}_\delta(x,y)+f_y(x))\nonumber\\
&=& {S}_\delta(\mathfrak{n},x)\: {S}_\delta(x,y) + f_{\mathfrak{n}}(x)^\ast \:{S}_\delta(x,y)
+ {S}_\delta(x,\mathfrak{n})^\ast \:f_y(x) + f_{\mathfrak{n}}(x)^\ast\: f_y(x)\:,
\end{eqnarray}
where again $y\in B_{\delta/2}(\mathfrak{n})$ and $x\in\tilde{M}$. Let us
analyze the $x$-integrals of the obtained expressions.
The integral over the first term gives precisely~$G_\delta$, (\ref{Gddef}).
For the last three expressions we can interchange the integral
with the limit $y\to \mathfrak{n}$, as the next Lemma shows.

\begin{Lemma} \label{2}
The following limits can be taken inside the integral,
\begin{eqnarray}
\lim_{\mathfrak{n} \neq y\to \mathfrak{n}}
\int_{\tilde{M}}f_{\mathfrak{n}}(x)^\ast\: {S}_\delta(x,y)\: dx &=&
\int_{\tilde{M}} f_{\mathfrak{n}}(x)^\ast{S}_\delta(x,\mathfrak{n}) \:dx \label{b} \label{2.6} \\
\lim_{\mathfrak{n} \neq  y\to \mathfrak{n}}\int_{\tilde{M}}{S}_\delta(x,\mathfrak{n})^\ast\: f_y(x) \:dx
&=& \int_{\tilde{M}} {S}_\delta(x,\mathfrak{n}) \:f_{\mathfrak{n}}(x)\:dx \label{c}\label{2.7}\\
\lim_{\mathfrak{n} \neq y\to \mathfrak{n}} \int_{\tilde{M}} f_{\mathfrak{n}}(x)^\ast \:f_y(x)^\ast \:dx
&=& \int_{\tilde{M}}f_{\mathfrak{n}}(x) \:f_{\mathfrak{n}}(x)^\ast \:dx \:. \label{d}
\end{eqnarray}
\end{Lemma}
\Proof The equations (\ref{c}) and (\ref{d}) follow immediately
from Lebesgue's dominated convergence theorem using Lemma~ \ref{lemma1}
and the fact that the pole of the Green's function on the sphere is
integrable (see Section~\ref{secGFS} for details).
The proof of (\ref{b}) is a bit
harder, and we use the symmetry of $\ssn$ on $S^n_\sigma$: Let
$\varphi_y$ be an isometry on $S^n_\sigma$ with $\varphi_y(y)=\mathfrak{n}$.
Then, since the Lebesgue integral over $S^n_\sigma$ is invariant under
$\varphi$,
\[ \int_{\tilde{M}}f_{\mathfrak{n}}(x)^\ast{S}_\delta(\tilde{x},\tilde{y})d\tilde{x}
\;=\; \int_{S^n_\sigma} f_{\mathfrak{n}}(x)^\ast{S}_\delta(x,y)dx \;=\;
\int_{S^n_\sigma} f_{\varphi_y(\mathfrak{n})}(x)^\ast {S}_\delta(x,\mathfrak{n})\: dx \:. \]
Now we can again apply Lebesgue's dominated convergence theorem,
\begin{eqnarray*}
\lim_{\mathfrak{n} \neq y\to \mathfrak{n}}\int_{\tilde{M}}f_{\mathfrak{n}}(x)^\ast{S}_\delta(x,y)\:dx
&=& \lim_{\mathfrak{n} \neq y\to \mathfrak{n}} \int_{S^n_\sigma} f_{\varphi_y(\mathfrak{n})}(x)^\ast{S}_\delta(x,\mathfrak{n})\:
dx \\
&=& \int_{S^n_\sigma}f_{\mathfrak{n}}(x)^\ast{S}_\delta(x,\mathfrak{n})\:dx \;=\; \int_{\tilde{M}}
  f_{\mathfrak{n}}(x)^\ast{S}_\delta(x,\mathfrak{n})\:dx\: .
\end{eqnarray*}

\vspace*{-1.85em} \QED

We write~(\ref{h}) in the form
\begin{eqnarray*}
\lefteqn{ \tilde{S}_{\bar{M}}(\mathfrak{n},x)\: \tilde{S}_{\bar{M}}(x,y)-{S}_\delta(\mathfrak{n},x)\: {S}_\delta(x,y) } \\
&=& F_y(x) \;:=\; f_{\mathfrak{n}}(x)^\ast\: {S}_\delta(x,y) + {S}_\delta(x,\mathfrak{n})^\ast\:
f_y(x) +f_{\mathfrak{n}}(x)^\ast \:f_y(x) \:.
\end{eqnarray*}
According to Lemma~\ref{lemma1} and Lemma~\ref{2},
we may commute the integral over
$x$ with the limit $y\to \mathfrak{n}$ as follows,
\begin{eqnarray*}
\lefteqn{ \int_{\tilde{M}}(\tilde{S}_{\bar{M}}(\mathfrak{n},x)\: \tilde{S}_{\bar{M}}(x,\mathfrak{n})
-{S}_\delta(\mathfrak{n},x){S}_\delta(x,\mathfrak{n})) \: dx } \\
&=& \int_{\tilde{M}} \lim_{y\to \mathfrak{n}} F_y(x)\:
dx \;=\; \lim_{y\to \mathfrak{n}} \int_{\tilde{M}} F_y(x)\:dx \\
&=&\lim_{y\to \mathfrak{n}} \int_{\tilde{M}}
(\tilde{S}_{\bar{M}}(\mathfrak{n},x)\tilde{S}_{\bar{M}}(x,y)-{S}_\delta(\mathfrak{n},x){S}_\delta(x,y))\: dx \:.
\end{eqnarray*}
This concludes the proof of Theorem~\ref{theorem1}.

\section{Pointwise Estimate of $G$ near the Pole} \label{sec5}
In this section we shall estimate the quantity~$\lim_{y \to \mathfrak{n}}
(G(\mathfrak{n},y)-G_{\delta}(y))$ appearing in Theorem~\ref{theorem1}.
The first difficulty is that that~$G_\delta$ is defined by
an integral~(\ref{Gddef}) and is thus rather complicated;
it would be more convenient to work instead with the
Green's function on the sphere~$G_{S^n_\sigma}(\mathfrak{n},y)$.
Therefore, we introduce the function~$H_\delta$ by
\begin{equation} \label{Hddef}
H_\delta(y) \;=\; G_{S^n_\sigma}(\mathfrak{n},y) - G_\delta(y) \:.
\end{equation}
This function depends only on the Green's functions on the sphere
and can be computed explicitly (for details
see Lemma~\ref{lemma3} below).
Thus it remains to control the difference of
the Green's functions~$G$ and~$G_{S^n_\sigma}$.
First we need to localize the last Green's function inside the
spherical cap, so that we can lift it to~$\bar{M}$.
To this end, we multiply it with the function~$\eta$ (\ref{edef}),
which is supported inside the spherical cap and is identically
equal to one in a neighborhood of the north pole.
Then our task is to estimate the limit
\begin{equation} \label{gdef}
\lim_{y \to \mathfrak{n}} \gamma(y) \quad {\mbox{with}} \quad
\gamma(y) \;:=\; G(\mathfrak{n},y) - \eta(y)\:
G_{S^n_\sigma}(\mathfrak{n},y)\:.
\end{equation}

Our strategy is as follows. When we apply the Dirac operator
squared to~$\gamma$, the $\delta$-contributions cancel,
\begin{equation} \label{hdef}
h(y) \;:=\; \tilde{\D}^2 \gamma(y) \;=\;
\tilde{\D}^2 \left( \eta(y)\: G_{S^n_\sigma}(\mathfrak{n},y) \right)
- \eta(y) \left(\tilde{\D}^2  G_{S^n_\sigma}(\mathfrak{n},y) \right) .
\end{equation}
The resulting terms all involve derivatives of~$\eta$
and are thus supported in the annulus~$\frac{\delta}{2}
< d <\delta$, where~$G_{S^n_\sigma}$ is smooth (for details see
Section~\ref{secGFS}). We conclude that~$h \in C^\infty_0(C)$.
We consider the equation
\begin{equation} \label{ell1}
\tilde{\D}^2 \gamma \;=\; h \qquad {\mbox{on $\bar{M}$}}
\end{equation}
as an elliptic equation for~$\gamma$. Standard elliptic regularity
theory yields that~$\gamma \in C^\infty(\bar{M})$.
Furthermore, the operator~$\tilde{\D}^2$ is essentially self-adjoint
on the Hilbert space~$L^2(\bar{M})$ of square-integrable
spinors with domain~$C^\infty(\bar{M})$.
According to the Lichnerowicz-Weitzenb{\"o}ck formula and the
fact that the scalar curvature is non-negative on~$\bar{M}$,
we know that the operator~$\tilde{\D}^2$ is strictly positive,
and we obtain
\begin{equation} \label{infest}
\|\gamma\|_{L^2(\bar{M})} \;\leq\; \frac{1}{\inf \spec(\tilde{\D}^2)}\:
\|h\|_{L^2(\bar{M})} \:.
\end{equation}
This~$L^2$-estimate is clearly not good enough, we need
a {\em{pointwise}} estimate. The general method is to derive
integral estimates for the derivatives of~$\gamma$ and then to
apply the Sobolev imbedding theorem.
The Sobolev imbedding theorem on a manifold involves
the isoperimetric constant (see e.g.\ \cite{He, FK}).
The basic problem is that the isoperimetric constant on~$\bar{M}$
depends on the unknown geometry in the compact set~$K$
and is therefore not under control.
In order to bypass this problem, we shall always work with
functions which are supported inside the spherical cap,
so that we can use the Sobolev imbedding on~$S^n_\sigma$.
More specifically, we work with the Sobolev inequality~\cite{He}
\begin{equation} \label{sobolev}
|\gamma(\mathfrak{n})| \;\leq\; \sup_{S_\sigma^n} |\eta^k \gamma| \;\leq\;
c_S\: \|\eta^k \gamma\|_{H^{k,2}(S^n_\sigma)} \quad{\mbox{for}}\quad
k \;>\; \frac{n}{2}\:.
\end{equation}
Here~$c_S$ is the Sobolev constant on the unit sphere~$S^n$, and
the Sobolev norm~$\|.\|_{H^{k,2}(S^n_\sigma)}$ is defined by
\[ \|f\|^2_{H^{k,2}(S^n_\sigma)} \;=\; \sum_{\kappa {\mbox{\scriptsize{ with }}} |\kappa|\leq k} \sigma^{2 |\kappa| - n}
\int_{S^n_\sigma} \|\nabla^\kappa f(x)\|^2 \:dx \:. \]
We inserted the factor~$\sigma^{2 |\kappa| - n}$ for convenience; it makes the
Sobolev norm invariant under scalings of~$\sigma$.

It remains to get estimates for the Sobolev norms in~(\ref{sobolev}).
The next lemma shows that we can equivalently
consider the $L^2$-norms of higher powers of Dirac operator.
\begin{Lemma} \label{lemma61}
There is a constant~$c$ which depends only on~$n$
and the quotient~$\delta/ \sigma$ (but is independent of~$\gamma$) such that
\begin{equation} \label{e1}
\|\eta^{k+1} \gamma\|^2_{H^{k,2}(S^n_\sigma)} \;\leq\; c \sum_{l=0}^k \sigma^{2 l - n} \;\|\eta^{l+1}\: \tilde{\D}^l \gamma\|_{L^2(S^n_\sigma)}^2 \:.
\end{equation}
\end{Lemma}
{\Proof} Since both sides of the inequality have the same
scaling in~$\sigma$, we may assume that~$\sigma=1$.
Using the Leibniz rule and the boundedness
of~$\eta$ and its derivatives, we immediately get
\begin{equation} \label{e2}
\|\eta^{k+1} \gamma\|^2_{H^{k,2}(S^n)} \;\leq\; c \sum_{|\kappa|\leq k}
\int_{S^n_\sigma} \|\nabla^\kappa (\eta^{|\kappa|+1} \gamma(x))\|^2 \:dx\:.
\end{equation}
Thus it suffices to show that the right side of~(\ref{e2}) is
controlled by the right side of~(\ref{e1}). We proceed by induction in~$k$.
For~$k=0$ there is nothing to prove. Suppose that the inequality holds
for given~$k$. Then, by the Leibniz rule and the Schwarz inequality,
\[ \|\eta^{k+2}\: \tilde{\D}^{k+1} \gamma\|_{L^2(S^n)}^2 \;=\; 
\|\tilde{\D}^{k+1} (\eta^{k+2} \gamma)\|_{L^2(S^n)}^2 + {\mbox{(l.o.t)}} \:, \]
where ``(l.o.t.)'' stands for $L^2$-norms of $\nabla^{\kappa}(\eta^{|\kappa|} \gamma)$ for
lower orders $|\kappa|<k+1$. Using the induction hypothesis, we obtain
\[ \sum_{l=0}^{k+1} \|\eta^{l+1}\: \tilde{\D}^l \gamma\|_{L^2(S^n)}^2 \;\geq\; \frac{1}{c}\: \sum_{l=0}^{k+1} 
\|\tilde{\D}^l (\eta^{l+1} \gamma)\|_{L^2(S^n)}^2\:. \]
The Lichnerowicz-Weitzenb{\"o}ck formula together with the fact that the scalar
curvature is non-negative imply that~$\tilde{D}^2 \geq \Delta$. Hence,
setting~$\chi = \eta^{k+1}$,
\[ \langle\tilde{\D}^{k+1}\: \chi \gamma, \tilde{\D}^{k+1}\: \chi \gamma \rangle_{L^2(S^n)}
\;=\;  \langle\tilde{\D}^{2k+2}\: \chi \gamma, \chi \gamma \rangle_{L^2(S^n)}
\;\geq\; \langle\Delta^{k+1}\: \chi \gamma, \chi \gamma \rangle_{L^2(S^n)} \:. \]
Integrating by parts, of each Laplacian we bring one derivative
on each side of the scalar product. Commuting covariant derivatives
we pick up curvature terms, which are clearly bounded on~$S^n$.
We thus obtain
\[ \langle\Delta^{k+1}\: \chi \gamma, \chi \gamma \rangle_{L^2(S^n)} \;=\;
\sum_{|\kappa|=k+1} \langle\nabla^{\kappa}\: \chi \gamma, \nabla^{\kappa} \chi \gamma \rangle_{L^2(S^n)} + {\mbox{(l.o.t.)}} \:. \]
Again using the induction hypothesis, the result follows.
\QED

The~$L^2$-norms on the right side of~(\ref{e1}) can easily be
estimated similar to~(\ref{infest}),
\begin{eqnarray*}
\|\eta^{l+1}\: \tilde{\D}^l \gamma\|_{L^2(S^n_\sigma)}^2 &\leq&
\|\tilde{\D}^l \gamma\|_{L^2(\bar{M})}^2 \\
&\leq& \frac{1}{\inf \spec(\tilde{\D}^2)^2}\: \|\tilde{\D}^{l+2} \gamma\|_{L^2(\bar{M})}^2
\;=\; \frac{1}{\inf \spec(\tilde{\D}^2)^2}\: \|\tilde{\D}^{l} h\|_{L^2(S^n_\sigma)}^2 \:.
\end{eqnarray*}
Notice that the norm on the very right depends only on the geometry
of the spherical cap.
Combing this last estimate with~(\ref{sobolev}) and Lemma~\ref{lemma61},
we obtain the following result.
\begin{Corollary} \label{cor52}
There is a constant~$c$ depending only on~$n$ and $\delta/ \sigma$ such that
\[ \lim_{y \to \mathfrak{n}} \gamma(y) \;\leq\; \frac{c\: \sigma^{-n}}
{\inf \spec(\tilde{\D}^2)}\:. \]
\end{Corollary}
{\Proof} Putting together the above estimates, we obtain
\begin{eqnarray}
\lim_{y \to \mathfrak{n}} \gamma(y)^2 &=& \gamma(\mathfrak{n})^2 \;\leq\;
c \sum_{l=0}^k \sigma^{2 l - n} \;\|\eta^{l+1}\: \tilde{\D}^l \gamma\|_{L^2(S^n_\sigma)}^2 \nonumber \\
&\leq& \frac{c\: \sigma^{-2n}}{\inf \spec(\tilde{\D}^2)^2}
\sum_{l=1}^k \sigma^{2 l+n} \;\|\tilde{\D}^{l} h\|_{L^2(S^n_\sigma)}^2 \:.
\label{eqn}
\end{eqnarray}
Since~$h$ is a given function on the spherical cap,
the summands in the last sum can be bounded by a
constant depending only on~$n$, $\delta$ and $\sigma$.
In order to determine the scaling in~$\sigma$, we first
note that
\[ \tilde{\D}^2_x G_{S^n_\sigma}(x,y) \;=\; \delta_{S^n_{\sigma}}(x,y) \;=\;
\sigma^{-n} \:\delta_{S^n_1}(\sigma^{-1} x, \sigma^{-1} y) \;=\;
\sigma^{-n}\:\tilde{\D}^2_x G_{S^n_1}(\sigma^{-1} x, \sigma^{-1} y) \:, \]
and thus
\[ G_{S^n_\sigma}(x,y) \;=\; \sigma^{-n+2}\: G_{S^n_1}(\sigma^{-1} x, \sigma^{-1} y)\:. \]
Using this in~(\ref{hdef}), one sees that~$h$ scales
like~$\sigma^{-n}$. We conclude that the last sum in~(\ref{eqn})
is scaling invariant, and therefore
this sum can be bounded by a constant which depends only
on the quotient $\delta/ \sigma$.
\QED

\section{The Green's Functions on the Sphere} \label{secGFS}
In the previous constructions we used the Green's functions
on the sphere~$S_{S^n_\sigma}$ and $G_{S^n_\sigma}$. We shall now
compute these Green's functions and estimate the composite
expressions~$S_\delta(x,\mathfrak{n}) \:S_\delta(\mathfrak{n}, x)$
and~$G_\delta - G_{S^n_\sigma}$.

Under the conformal transformation of $S^n_\sigma \setminus \{\mathfrak{n}\}$ to
Euclidean~$\R^n$, the Green's function~$S_{S^n_\sigma}$ clearly goes
over to the Green's function of Euclidean space~(\ref{SEuclid}).
Applying Lemma~\ref{conformal change}, we thus obtain for the
Green's function on~$S^n_\sigma$ in the coordinates of the stereographic
projection from the north pole the explicit formula
\begin{equation} \label{SSN}
S_{S^n_\sigma}(x,y) \;=\; -\frac{1}{\omega_{n-1}}\: \frac{x-y}{|x-y|^n}\:
\left( \frac{2 \sigma^2}{\sigma^2+|x|^2} \right)^{\frac{1-n}{2}}
\left( \frac{2 \sigma^2}{\sigma^2+|y|^2} \right)^{\frac{1-n}{2}} .
\end{equation}
In particular, one sees that the Green's function on the sphere is
smooth away from the diagonal, and that the pole at~$x=y$ is integrable.

In the next lemma we compute the product~${S}_\delta(x,\mathfrak{n}) {S}_\delta(\mathfrak{n},x)$
as well as~$G_\delta$ and $H_\delta$ defined by~(\ref{Gddef}, \ref{Hddef}).
It is most convenient to work on~$S^n_\sigma$
in the coordinate system obtained by stereographic projection from
the south pole. The corresponding radial coordinate~$r'$ is related
to the radial coordinate~$r$ in the stereographic projection from the north
pole by
\begin{equation} \label{rpdef}
r' \;=\; \frac{\sigma^2}{r}\:.
\end{equation}

\begin{Lemma}\label{lemma3}
Setting~$R'=\sigma^2/R$, the following identities hold inside
the spherical cap~$C$,
\begin{eqnarray}
{S}_\delta(x,\mathfrak{n})\: {S}_\delta(\mathfrak{n},x) &=&
\frac{(2r'^2)^{1-n}}{\omega_{n-1}^2} \left(
\frac{2 \sigma^2}{\sigma^2 + r'^2} \right)^{1-n}\: \1 \label{Sdrel} \\
G_\delta(x) &=& \frac{1}{\omega_{n-1}}
\left(\frac{4 \sigma^2}{\sigma^2+r'^2}\right)^{\frac{1-n}{2}}
\int^{R'}_{r'} \frac{2\sigma^2}{\sigma^2 + \tau^2}
\:\frac{1}{\tau^{n-1}} \: d\tau \label{Gdrel} \\
H_\delta(x) &=& \frac{1}{\omega_{n-1}} \left(\frac{4 \sigma^2}
{\sigma^2+r'^2}\right)^{\frac{1-n}{2}} \int_{R'}^\infty
\frac{2\sigma^2}{\sigma^2 + \tau^2} \:\frac{1}{\tau^{n-1}} \: d\tau \:.
\label{Hdrel}
\end{eqnarray}
\end{Lemma}
{\Proof} The identity~(\ref{Sdrel}) can be obtained in two ways.
Either one computes the product~$S_\delta(x,y)\: S_\delta(y,x)$
with~$S_\delta$ in the ``north pole chart''~(\ref{SSN}) and takes the limit
$y \to \infty$. Alternatively, one can compute the product $S_\delta(x,\mathfrak{n})\: S_\delta(\mathfrak{n},x)$ in the chart of the stereographic
projection from the south pole.

For the rest of the proof we work with the stereographic
projection from the south pole, which we denote by
$\pi :S^n_c\to\mathbb{R}^n$.
Using the explicit formulas~(\ref{SEuclid}, \ref{SSN}, \ref{gSn}), we
obtain
\begin{eqnarray*}
G_\delta(y) &=&
\int_{S^n_\sigma} \:{S}_\delta(\mathfrak{n},x){S}_\delta(x,y) \:dx \;=\;
\int_{B_\delta(\mathfrak{n})} S_{S_\sigma^n}(\mathfrak{n},x)\: S_{S_\sigma^n}(x,y)\: dx \\
&=& \int_{\pi(B_\delta(\mathfrak{n}))} S_{\mathbb{R}^n}(0,x) \:S_{\mathbb{R}^n}(x,y)\;
\frac{2 \sigma^2}{\sigma^2+|x|^2} \left(\frac{4 \sigma^2}{\sigma^2 + |y|^2}\right)^{\frac{1-n}{2}} dx \\
&=& \left(\frac{4 \sigma^2} {\sigma^2+|y|^2}\right)^{\frac{1-n}{2}}
\int_{\pi(B_\delta(\mathfrak{n}))} \left(\D_{\mathbb{R}^n,x} F(x) \right)
S_{\mathbb{R}^n}(x,y) \:dx
\end{eqnarray*}
with
\begin{equation} \label{Fdef}
F(x) \;:=\; \frac{1}{\omega_{n-1}} \int^{R'}_{r'(x)}
\frac{2\sigma^2}{\sigma^2 + \tau^2} \:\frac{1}{\tau^{n-1}} \: d\tau \:.
\end{equation}
We now integrate by parts. The boundary terms drop out because~$F$
vanishes on $\partial B_\delta(\mathfrak{n})$ (note that the pole at
the origin is of order~$(r')^{n-2}$, and so we do not need to worry about
boundary terms there). We conclude that
\[ G_\delta(y) \;=\; \left(\frac{4 \sigma^2}{\sigma^2+|y|^2}
\right)^{\frac{1-n}{2}} \int_{B_\delta(\mathfrak{n})} F(x)
\: \D_{\R^n_x} S_{\mathbb{R}^n}(x,y)\: dx
\;=\; \left(\frac{2 \sigma^2}{1+|y|^2}\right)^{\frac{1-n}{2}}\: F(y) \:, \]
proving~(\ref{Gdrel}). $H_\delta$ can be computed similarly,
\begin{eqnarray*}
H_\delta(y) &=&
\int_{S^n_\sigma \setminus B_\delta(\mathfrak{n})}
S_{S_\sigma^n}(\mathfrak{n},x)\: S_{S_\sigma^n}(x,y)\: dx \\
&=& \left(\frac{4 \sigma^2} {\sigma^2+|y|^2}\right)^{\frac{1-n}{2}}
\int_{\R^n \setminus \pi(B_\delta(\mathfrak{n}))} \left(\D_{\mathbb{R}^n,x} F(x) \right)
S_{\mathbb{R}^n}(x,y) \:dx\:.
\end{eqnarray*}
Now we integrate by parts. Since~$\D_{\mathbb{R}^n,x} S_{\mathbb{R}^n}(x,y)$
gives a contribution only for~$x=y \not \in \R^n \setminus \pi(B_\delta(\mathfrak{n}))$,
only the boundary integrals contribute. The boundary terms on
$\partial (\pi (B_\delta(\mathfrak{n})))$ again drop out because~$F$ vanishes there.
Thus it remains to consider the boundary terms at infinity,
\[ H_\delta(y) \;=\; \left(\frac{4 \sigma^2}
{\sigma^2+|y|^2}\right)^{\frac{1-n}{2}} \:
\lim_{L \to \infty}
\int_{S^{n-1}_L} \nu \cdot F(x)\: S_{\mathbb{R}^n}(x,y) \:dx \:, \]
where~$\nu \cdot$ denotes Clifford multiplication with the outer normal.
Putting in~(\ref{SEuclid}, \ref{Fdef}), we obtain~(\ref{Hdrel}).
\QED

Substituting~(\ref{Gdrel}, \ref{Hdrel}) into~(\ref{Hddef}), we
also obtain an explicit expression for the Green's function~$G_{S^n_\sigma}$,
\[ G_{S^n_\sigma}(\mathfrak{n}, x) \;=\;  \frac{1}{\omega_{n-1}}
\left(\frac{4 \sigma^2}{\sigma^2+r'^2}\right)^{\frac{1-n}{2}}
\int^{\infty}_{r'} \frac{2\sigma^2}{\sigma^2 + \tau^2} \:\frac{1}{\tau^{n-1}} \: d\tau\:. \]
This formula shows that~$G_{S^n_\sigma}(x,y)$ is indeed smooth away from
the diagonal. 

\section{A Weighted~$L^1$-Estimate of the Deviation Operator} \label{sec7}
In this section we first combine the previous results to
derive an integral estimate for the trace of the spinor
operator (Theorem~\ref{Thmw1}). We then introduce the
so-called deviation operator, which gives us information on
how much the spinor operator differs in the asymptotic end from the
spinor operator in Schwarzschild. Using a positivity argument
(Lemma~\ref{partialwave}), we can then prove the main results of this paper:
a weighted~$L^1$-estimate of the deviation operator (Theorem~\ref{thmspop})
and a weighted~$L^2$-estimate of the Witten spinors (Theorem~\ref{thmwitten}).

Let~$\mu$ be a function which in the asymptotic end coincides with
the norm of the spinor operator in Schwarzschild~(\ref{PiSch}) and
vanishes otherwise,
\begin{equation} \label{mudef}
\mu \;:\; M \to \R\:,\quad
\mu(x) \;=\; \chi_{M \setminus K}(x) \:\left(1 + \frac{1}{r(x)^{n-2}}
\right)^{-2\: \frac{n-1}{n-2}} \:.
\end{equation}
In the next Theorem we compute an integral involving
the trace of the matrix~$\Pi-\mu \1$.
\begin{Thm} \label{Thmw1}
The spinor operator on~$M$ satisfies the following identity,
\begin{equation} \label{eqw1}
\int_M \Tr \left( \Pi(x) - \mu(x) \:\1 \right) \:
\lambda(x)\: dx \;=\; \omega^2_{n-1}\:(2 \sigma^2)^{n-1}
\lim_{y \to \mathfrak{n}} \Tr \left( \gamma(y) \right)
+ N \alpha \:,
\end{equation}
with~$\gamma$ according to~(\ref{gdef}), and where~$\alpha$
is given by explicit integrals in Euclidean space,
\begin{equation} \label{adef}
\alpha \;=\; \int_{B_R(0)} \frac{2 \sigma^2}{\sigma^2+|x|^2}\: d^nx
\:-\: \int_{B_R(0) \setminus \overline{B_\rho(0)}}
\left( 1 + \frac{1}{|x|^{n-2}}
\right)^{\frac{2}{n-2}} \:\lambda(\phi^{-1}(x))\: d^nx\:.
\end{equation}
Here~$d^nx$ is the Lebesgue measure in~$\R^n$, and~$\lambda$
is to be chosen as in Lemma~\ref{lemmacap}.
\end{Thm}
We point out that the parameters~$\gamma$ and~$\alpha$
clearly depend on the conformal compactification, but they are
(for a given function~$\lambda$) independent of~$R$.
This is obvious for~$\gamma$ because its definition~(\ref{gdef})
involves only the Green's functions on~$\bar{M}$ and~$S^n_\sigma$.
For the parameter~$\alpha$ it follows from the fact that for~$r>R$,
$\lambda$ is given by~(\ref{ldef}), so that the second integrand
in~(\ref{adef}) reduces to the first. Hence the integrals in~(\ref{adef})
remain unchanged if~$R$ is increased. \\[.5em]
{\em{Proof of Theorem~\ref{Thmw1}.}} Using Definition~\ref{defctso} and
Theorem~\ref{1.Theorem}, we get
\begin{eqnarray*}
\lefteqn{ \int_M \Tr \left( \Pi(x) - \mu(x) \:\1 \right) \:
\lambda(x)\: dx
\;=\; \int_{\ti{M}} \Tr \left( \tilde{\Pi}(x)  -
\mu(x)\: \lambda^{1-n}(x) \1 \right) \: dx } \\
&=& \int_{\ti{M}} \Tr \left( \omega^2_{n-1}\:(2 \sigma^2)^{n-1} \:\ti{S}_{\bar{M}}(x,\mathfrak{n}) \:\ti{S}_{\bar{M}}(\mathfrak{n},x)
\:-\: \mu(x)\: \lambda^{1-n}(x)\:\1 \right) \: dx \\
&=&  \omega^2_{n-1}\:(2 \sigma^2)^{n-1} \int_{\ti{M}}
\Tr \left( \ti{S}_{\bar{M}}(\mathfrak{n},x) \:\ti{S}_{\bar{M}}(x,\mathfrak{n}) 
\:-\: S_\delta(\mathfrak{n},x)\: S_\delta(x,\mathfrak{n}) \right) dx \\
&&+ \int_{\ti{M}} \Tr \left( \omega^2_{n-1}\:
(2 \sigma^2)^{n-1}\: S_\delta(x, \mathfrak{n})\: S_\delta(\mathfrak{n}, x)
\:-\: \mu(x)\: \lambda^{1-n}(x)\:\1 \right) dx \:.
\end{eqnarray*}
According to~(\ref{ldef}, \ref{mudef})
and~(\ref{Sdrel}, \ref{rpdef}),
\[ \omega^2_{n-1}\:(2 \sigma^2)^{n-1}\: S_\delta(x, \mathfrak{n})\: S_\delta(\mathfrak{n}, x) \:-\: \mu(x)\:\lambda^{1-n}(x)
 \;\equiv\; 0 \qquad {\mbox{on $C$}} \:, \]
and thus the last integral reduces to an integral over the
annular region $\bar{M} \setminus (K \cup C)$. 
Furthermore, we can apply Theorem~\ref{theorem1} as well as~(\ref{Hddef}, \ref{gdef}) to obtain
\begin{eqnarray*}
\lefteqn{ \int_M \Tr \left( \Pi(x) - \mu(x) \:\1 \right) \:
\lambda(x)\: dx
\;=\; \omega^2_{n-1}\:(2 \sigma^2)^{n-1}
\lim_{y \to \mathfrak{n}} \Tr \left( \gamma(y) \right) } \\
&&+ \omega^2_{n-1}\:(2 \sigma^2)^{n-1}\: \Tr (H_\delta(\mathfrak{n}))
\;-\; \int_{M \setminus (K \cup C)} \Tr \left( 
\mu(x)\: \lambda(x)\:\1 \right) dx \:.
\end{eqnarray*}
The last two terms can be computed explicitly with~(\ref{Hdrel})
and~(\ref{Schmet}),
\begin{eqnarray*}
\omega^2_{n-1}\:(2 \sigma^2)^{n-1}\: H_\delta(\mathfrak{n})
&=& \omega_{n-1}\: \sigma^{2(n-1)}
\int_{R'}^\infty \frac{2 \sigma^2}{\sigma^2 + \tau^2}\:\tau^{n-1}\: d\tau \\
&=& \omega_{n-1} \int_0^R\: \frac{2 \sigma^2}{\sigma^2+r^2}\: r^{n-1}\:dr
\;=\; \int_{B_R(0)} \frac{2 \sigma^2}{\sigma^2+|x|^2}\: d^nx \\
\int_{M \setminus (K \cup C)} \mu(x)\: \lambda(x)\: dx
&=& \int_{B_R(0) \setminus \overline{B_\rho(0)}}
\left( 1 + \frac{1}{|x|^{n-2}} \right)^{\frac{2}{n-2}}
\:\lambda(\phi^{-1}(x))\: d^nx\:.
\end{eqnarray*}

\vspace*{-1.65em} \QED

We now analyze the integral in~(\ref{eqw1}) in more detail, with
the aim of getting a connection to an~$L^1$-norm.
Inside the compact set~$K$, the function~$\lambda$ vanishes.
Thus according to Definition~\ref{defso}, for every
spinor~$\psi \in \Sigma_x$,
\begin{equation} \label{cant}
\langle \psi,  (\Pi(x) - \mu(x) \:\1)\: \psi \rangle
\;=\; \langle \psi,  \Pi(x) \: \psi \rangle \;=\;
\sum_{i=1}^N |\langle \psi_i(x), \psi \rangle|^2 \;\geq\; 0 \:.
\end{equation}
Hence the matrix~$\Pi - \mu \1$ is positive, and we can
control the sup-norm by the trace,
\begin{equation} \label{suptrace}
\| \Pi(x) - \mu(x) \:\1 \| \;\leq\; \Tr \left(
\Pi(x) - \mu(x) \:\1 \right) \spc {\mbox{for all $x \in K$}}.
\end{equation}
In the asymptotic end, where~$\mu>0$, we
cannot expect that the operator~$\Pi - \mu \1$ is still positive.
Nevertheless, the next lemma shows that the integral over the
trace is indeed positive and can be identified with the trace
of a positive operator, which we call deviation operator.

\begin{Def} \label{defdo}
Working in the asymptotic end in the chart~$(\phi,\, M \setminus K)$,
we introduce for every solution~$\psi_i$ of the boundary value problem~(\ref{bv}) the {\bf{Witten deviation}}~$\delta \psi_i$ by
\[ \delta \psi_i \;=\; \psi_i - 
\left( 1 + \frac{1}{r(x)^{n-2}}
\right)^{-\frac{n-1}{n-2}} \psi_{0,i} \:. \]
For every~$x \in M \setminus K$, the {\bf{deviation operator}}~$\delta \Pi$
is defined by
\[ \delta \Pi(x) \::\: \Sigma_x \to \Sigma_x \::\:
\psi \;\mapsto\; \sum_{i=1}^N \langle \delta \psi_i(x), \psi \rangle
\:\delta \psi_i(x) \:. \]
\end{Def}
Thus the deviation operator is defined similar to the spinor operator;
one only replaces the Witten spinors by the corresponding
Witten deviations. Repeating the argument in~(\ref{cant}),
one sees that the deviation operator is also positive.

\begin{Lemma} \label{partialwave}
\[ \int_{M \setminus K} \Tr \left( \Pi(x) - \mu(x) \:\1 \right) \:
\lambda(x)\: dx
\;=\; \int_{M \setminus K} \Tr \left( \delta \Pi \right) \:
\lambda(x)\: dx \:. \]
\end{Lemma}
{\Proof}
We work in the chart~$(\phi,\, M \setminus K)$
and choose in~$\phi(M \setminus K) = \R^n \setminus
\overline{B_{\rho}(0})$ polar coordinates~$(r, \omega)$ with~$\omega
\in S^{n-1}$. 
According to the behavior of the Dirac operator under conformal
transformations~(\ref{ctrans}),
every Witten spinor~$\psi$ can be written in the form
\[ \psi(x) \;=\; \left( 1 + \frac{1}{r(x)^{n-2}}
\right)^{-\frac{n-1}{n-2}}\: \Psi(x) \]
with~$\Psi$ a harmonic spinor on~$\R^n \setminus
\overline{B_{\rho}(0)}$ endowed with the Euclidean metric.
We now expand~$\Psi$ in partial waves,
\[ \Psi(r, \omega) \;=\; \Psi_0 + \sum_{l=1}^\infty
\Psi^l(\omega)\: \frac{1}{r^{l+n-2}} \:, \]
where~$l$ are the angular quantum numbers, and
the~$\Psi^l(\omega)$ are linear combinations of the corresponding
spin-weighted spherical harmonics. From the
smoothness of~$\psi$ and the asymptotics at infinity,
it is clear that the sum converges in~$L^2\! \left(\R^n \setminus
\overline{B_\rho(0)} \right)^N$.
Consequently, the Witten spinors and the Witten deviations have the following
partial wave expansions,
\begin{eqnarray*}
\psi_i(r, \omega) &=& \left( 1 + \frac{1}{r(x)^{n-2}}
\right)^{-\frac{n-1}{n-2}} \left(\psi_{0,i} + \sum_{l=1}^\infty
\Psi_i^l(\omega)\: \frac{1}{r^{l+n-2}} \right) \\
\delta \psi_i(r, \omega) &=& \left( 1 + \frac{1}{r(x)^{n-2}}
\right)^{-\frac{n-1}{n-2}} \:\sum_{l=1}^\infty
\Psi_i^l(\omega)\: \frac{1}{r^{l+n-2}} \:.
\end{eqnarray*}
Using that the spin-weighted spherical harmonics for different~$l$ are
orthogonal on~$L^2(S^{n-1})$, a short calculation shows that
for all~$r>R$,
\[ \int_{S^{n-1}} \Tr \left( \Pi(x) - \mu(x) \:\1 \right)(r, \omega)\: d\omega
\;=\; \int_{S^{n-1}} \Tr \left( \delta \Pi \right) (r, \omega)\: d\omega\:. \]

\vspace*{-1.65em} \QED

Combining Theorem~\ref{Thmw1} with the above lemma, we
immediately obtain the following identity for the
weighted~$L^1$-norm of the spinor operator and the deviation operator.
\begin{Corollary} \label{corend}
The spinor operator on~$M$ satisfies the following identity,
\[ 
\int_K \Tr ( \Pi(x)) \: dx
\:+\: \int_{M \setminus K} \Tr ( \delta \Pi(x) ) \:
\lambda(x)\: dx \;=\; \omega^2_{n-1}\:(2 \sigma^2)^{n-1}
\lim_{y \to \mathfrak{n}} \Tr \left( \gamma(y) \right)
+ N \alpha \:, \]
where~$\gamma$ and~$\alpha$ are given by~(\ref{gdef}).
and~(\ref{adef}).
\end{Corollary}

Putting in the estimate of Corollary~\ref{cor52}
gives the following result.
\begin{Thm} \label{thmspop}
There is a constant~$c$ depending only on the dimension such that
\begin{equation}
\int_K \|\Pi(x)\|\: dx \:+\:
\int_{M \setminus K} \|\delta \Pi(x)\|\:\lambda(x)\:dx
\;\leq\; c\:\frac{(\rho+1)^n}{\sigma^2 \,\inf \spec({\tilde{\D}}^2)} \:,
\end{equation}
where~$\lambda$ is to be chosen according to Lemma~\ref{lemmacap}.
\end{Thm}
{\Proof} According to~Lemma~\ref{partialwave}
and the positivity of~$\Pi(x)$ and~$\delta \Pi(x)$, we know that
\[ \int_K \|\Pi(x)\|\: dx \:+\:
\int_{M \setminus K} \|\delta \Pi(x)\|\:\lambda(x)\:dx
\;\leq\; \int_M \Tr \left( \Pi(x) - \mu(x) \:\1 \right) \:
\lambda(x)\: dx\:. \]
We now apply Theorem~\ref{Thmw1}, Corollary~\ref{cor52} and use that, according
to~(\ref{scaling}), $\sigma$ scales like the radius~$\rho$.
Furthermore, it is obvious from~(\ref{adef}, \ref{scaling})
that~$\alpha$ scales like~$\rho^n$. This gives the estimate
\[ \int_K \|\Pi(x)\|\:dx \:+\:
\int_{M \setminus K} \|\delta \Pi(x)\|\:\lambda(x)\:dx
\;\leq\; \frac{c\: (\rho+1)^n}{\sigma^2\, \inf \spec({\tilde{\D}}^2)}
\:+\: c\: (\rho+1)^n\:. \]
We finally show that the second term on the right can be bounded
by the first. To this end, we need to bound the lowest eigenvalue
of~$\tilde{\D}^2$ from above: We choose a smooth wave function~$\psi$ which is
supported in the spherical cap and consider its Rayleigh quotient,
\[ \inf \spec({\tilde{\D}}^2) \;\leq\; \frac{\langle \tilde{\D} \psi, \tilde{\D} \psi \rangle_{L^2(\bar{M})}}
{\langle \psi, \psi \rangle_{L^2(\bar{M})}} \;\leq\; \frac{c}{\sigma^2}\:. \]

\vspace*{-1.65em} \QED

From this theorem one obtains weighted
$L^2$-estimates for all Witten spinors.
\begin{Corollary} \label{thmwitten}
There is a constant~$c$ depending only on the dimension such that
every Witten spinor~$\psi$ satisfies the weighted $L^2$-estimate
\[ \int_K \|\psi(x)\|^2 \:dx \:+\:
\int_{M \setminus K} \|\delta \psi(x)\|^2\:\lambda(x)\:dx
\;\leq\; c\: \frac{(\rho+1)^n}{\sigma^2\, \inf \spec({\tilde{\D}}^2)} \]
with~$\lambda$ according to Lemma~\ref{lemmacap}.
\end{Corollary}
{\Proof} We choose a basis of Witten spinors~$\psi_1,\ldots, \psi_n$
such that~$\psi_1 =\psi$. Then for all~$\phi \in \Sigma_x$,
\[ \langle \phi, \Pi(x)\:\phi \rangle \;=\;
\sum_{i=1}^N |\langle \psi_i(x), \phi \rangle|^2 \;\geq\;
|\langle \psi(x), \phi \rangle|^2 \:. \]
Taking the supremum over all unit spinors~$\phi$, we conclude
that~$\|\psi(x)\|^2 \leq \|\Pi(x)\|$.
In the same way, one sees that~$\|\delta \psi(x)\|^2 \leq \|\delta
\Pi(x)\|$.
\QED
Theorem~\ref{thmintro} and Theorem~\ref{thmintro2}
are a special case of Corollary~\ref{thmwitten}
and Corollary~\ref{corend}, respectively.\\[0em]

\noindent
{\em{Acknowledgments:}}
We would like to thank the ``Erwin-Schr{\"o}dinger-Institut,'' Wien,
and the ``Centre de Recherche
Math{\'e}mathiques,'' Montr{\'e}al, for hospitality and support.


\noindent
NWF I -- Mathematik,
Universit{\"a}t Regensburg, 93040 Regensburg, Germany, \\
{\tt{Felix.Finster@mathematik.uni-regensburg.de}}, \\
{\tt{Margarita.Kraus@mathematik.uni-regensburg.de}}

\end{document}